\documentclass[12pt]{amsart}

\usepackage{pb-diagram}

\usepackage{amsfonts}
\usepackage{amsmath}
\usepackage{amssymb}

\allowdisplaybreaks

\newcommand{\tr}{\textnormal{tr}}

\newcommand{\dbar}{\overline{\partial}}

\newcommand{\ddbar}{\sqrt{-1}\partial\dbar}

\newtheorem{theorem}{Theorem}[section]

\newtheorem{lemma}{Lemma}[section]

\newtheorem{definition}{Definition}[section]

\renewcommand{\thefootnote}{\fnsymbol{footnote}}
\newcommand{\starttext}{ \setcounter{footnote}{0}
\renewcommand{\thefootnote}{\arabic{footnote}}}

\newcommand{\beq}{\begin{equation}}
\newcommand{\bea}{\begin{eqnarray}}
\newcommand{\eea}{\end{eqnarray}} \newcommand{\ee}{\end{equation}}

\def\ba{\begin{eqnarray}}
\def\ea{\end{eqnarray}}

\def\ti\tilde

\def\u{\underline}

\def\tr{{\rm tr}}
\def\det{{\rm det}}

\def\ti{\tilde}

\newcommand{\bl}{\lambda}
\newcommand{\bmu}{\mu}

\begin{document}

\starttext \baselineskip=18pt \setcounter{footnote}{0}

\title[Uniform $L^\infty$ estimates: subsolutions]{Uniform $L^\infty$ estimates: subsolutions to fully nonlinear partial differential equations}

\author{Bin Guo and Duong H. Phong}

\thanks{Work supported in part by the National Science Foundation under grants DMS-2203273, and DMS-2303508, and the collaboration grant 946730 from Simons Foundation.}

\address{Department of Mathematics \& Computer Science, Rutgers University, Newark, NJ 07102}

\email{bguo@rutgers.edu}

\address{Department of Mathematics, Columbia University, New York, NY 10027}

\email{phong@math.columbia.edu}

\begin{abstract}

Uniform bounds are obtained using the auxiliary Monge-Amp\`ere equation method for solutions of very general classes of fully non-linear partial differential equations, assuming the existence of a ${\mathcal{C}}$-subsolution in the sense of G. Sz\'ekelyhidi and B. Guan. The main advantage over previous estimates is in the K\"ahler case, where the new estimates remain uniform as the background metric is allowed to degenerate.

{\footnotesize }

\end{abstract}

\maketitle

\baselineskip=15pt
\setcounter{equation}{0}
\setcounter{footnote}{0}

\section{Introduction}
\setcounter{equation}{0}

This paper is motivated by two relatively recent developments. The first is the new theory of ${\mathcal{C}}$-subsolutions for fully non-linear equations on compact Hermitian manifolds due to G. Sz\'ekelyhidi \cite{Sz} and B. Guan \cite{Gu}, a parabolic version of which was subsequently developed in \cite{PT}. Subsolutions are of great importance for non-linear equations, as their existence are clearly a necessary condition for the solvabiliy of the equations themselves.
A significant difficulty for equations on compact manifolds was that the naive notion may only allow trivial subsolutions, and this difficulty was only overcome recently with the notion of ${\mathcal{C}}$-subsolution introduced in \cite{Sz, Gu, PT}. The second development is the new auxiliary Monge-Amp\`ere equation method for $L^\infty$ estimates introduced by the authors in joint work with F. Tong for equations respectively on K\"ahler manifolds \cite{GPT} and Hermitian manifolds \cite{GP22a}. A great advantage of these new methods is that they lead to $L^\infty$ estimates which are uniform with respect to large classes of background metrics, the power of which can be appreciated from many recent applications in complex geometry \cite{Paun, GPS, GPSS22, GPSS23}.
A natural question is then whether the existence of subsolutions can lead, not just $L^\infty$ estimates, but $L^\infty$ estimates with uniformity comparable to those in \cite{GPS, GPSS22}. The goal of the present paper is to provide a positive answer to this question.

\medskip
We begin by providing a general formulation of the equations which we shall consider.
Let $X$ be a complex manifold of complex dimension $n$ with a Hermitian metric $\omega_X$.
Let $\omega$ be another Hermitian metric and  $\chi$ a closed $(1,1)$-form on $X$. For a smooth function $\varphi$, let $h_\varphi = \omega^{-1}\cdot \chi_\varphi: TX\to TX$ be the endomorphism defined by $\omega$ and $\chi_\varphi = \chi+\ddbar \varphi$, and let $\bl[h_\varphi]$ is the (unordered) vector of eigenvalues of this endomorphism. We consider the following fully nonlinear equation on $X$:
\bea\label{eqn:main}
f(\bl[h_\varphi]) = e^F, \quad {\mathrm{sup}}_X \varphi = 0,
\eea
where $e^F$ is a positive smooth function normalized such that $$\int_X e^{nF} \omega^n =[\omega]^n =  \int_X \omega^n. $$ 

The nonlinear operator $f: \Gamma\to {\mathbf R}_{>0}$ is required to  satisfy the following conditions: 

(i): $\Gamma\subset{\mathbf R}^n$ is an open convex symmetric cone such that $\Gamma_n\subset \Gamma \subset \Gamma_1$;


(ii): $\frac{\partial f}{\partial \lambda_i}(\bl) > 0$ for any $\bl = (\lambda_1,\ldots, \lambda_n)\in \Gamma$, which implies the ellipticity of the equation  (\ref{eqn:main}).

\medskip
In the papers \cite{GPT, AO}, it was further assumed that $f$ satisfies the structural condition \begin{equation}\label {eqn:f structure new}\prod_{i=1}^n\frac{\partial f}{\partial \lambda_i}(\bl) \ge \hat \gamma, \quad \forall \bl\in \Gamma_n\end{equation}
for some non-negative function $\hat \gamma$ which is required to be bounded from below by a fixed positive constant. This is a natural condition, which was subsequently shown by Harvey and Lawson \cite{HL22} to be satisfied by very large classes of equations, including all Garding-Dirichlet operators \cite{HL}. It is also a natural condition arising in the study of $C^2$ estimates \cite{AO, Dinew, CH}. However, as stressed in \cite{GPSS22}, in many important geometric applications, $\hat \gamma$ may vanish along a complex codimension one subvariety. Another important example of geometric equation which does not obey the condition of $\hat \gamma$ being strictly positive is the $J$ equation, which admits a solution if and only if a subsolution exists \cite{SW, We}. 

\medskip
We shall be interested in obtaining $L^\infty$ bounds for solutions of the equation (\ref{eqn:main}) assuming the existence of a subsolution, even in the absence of the structural condition (\ref{eqn:f structure new}) of the operator $f$. The following notion of $\mathcal C$-subsolution is taken from \cite{Sz} (see also \cite{PT}).
\begin{definition}\label{defn:2.1}
Given $\delta>0, R>0$, 
we say a $(1,1)$-form $\chi' = \chi - \ddbar \theta$ is a $(\delta,R)$-subsolution of the equation \eqref{eqn:main}, if the following holds
\begin{equation}\label{eqn:C sub}
(\lambda[h_{\chi'}] - \delta {\mathbf{1}} + \Gamma_n)\cap \partial\Gamma^{e^{F(z)}} \subset B(0, R)\subset {\mathbb R}^n,\quad \forall z\in X,
\end{equation}
where ${\mathbf 1} = (1,\ldots, 1)\in {\mathbb R}^n$ is the fixed vector,  $h_{\chi'} = \omega^{-1}\cdot \chi'$ is the endomorphism defined by $\chi'$ relative to $\omega$, and $$\Gamma^{e^{F(z)}} = \{ \lambda\in \Gamma~|~ f(\lambda)\ge e^{F(z)}  \}.$$ 

\end{definition}

The most important case for us is the K\"ahler case, to be discussed in Section 2, where estimates are obtained which are uniform with respect to large sets of forms $\chi$ reaching the boundary of the K\"ahler cone. The Hermitian case is discussed in Section 3. Finally, in Section 4, we discuss a notion of subsolution closer in spirit to numerical criteria in stabilty.

\section{The K\"ahler case} 
\label{section Kahler} 

In this section, we assume $(X,\omega_X)$ is a K\"ahler manifold. We will derive {\em uniform estimates} for $\varphi$ when a subsolution exists. Here  uniform estimates refer to the ones which do not depend on a fixed reference metric, in particular, it allows the reference metrics to degenerate. For the purpose of uniform estimates, we need to consider a subset of the space of K\"ahler metrics on $X$. As in the joint work \cite{GPSS22} with Song and Sturm,  for a given K\"ahler metric $\omega$, we will denote its relative volume function as
\begin{equation}\label{eqn:F omega}e^{F_\omega} = \frac{1}{V_\omega} \frac{\omega^n}{\omega_X^n},\quad V_\omega = \int_X \omega^n.\end{equation}
Given constants $A>0, K>0$, $p>n$ and a continuous nonnegative function $\gamma\in C^0(X)$, we consider the following subset of K\"ahler metrics
\bea
\label{eqn:Kset}
{\mathcal W}(A, K, p, \gamma, \omega_X) &: = & \Big\{ \omega \mbox { a K\"ahler metric} ~|~ [\omega]\cdot[\omega_X]^{n-1} \le A, e^{F_\omega}\ge \gamma, \\
&&\nonumber\quad \mbox{and }\, {\mathcal N}_p(\omega) = \int_X |F_\omega|^p e^{F_\omega}\omega_X^n \le K \Big\}.
\eea We will always assume that $\gamma$ satisfies ${\mathrm{dim}}_{{\mathcal H}} \gamma^{-1}(0) < 2n - 1$. 
When the parameters are fixed, we will simply write ${\mathcal W}$ as this set. 

Suppose $\chi' = \chi - \ddbar\theta$ is a $(\delta,R)$-subsolution to \eqref{eqn:main} as in Definition \ref{defn:2.1}, which satisfies the assumption:
\begin{equation}\label{eqn:chi}
- n\kappa_1 \le \tr_ {\omega} \chi' \le \kappa_2 , 
\end{equation}
for constants $\kappa_1\ge 0, \,\kappa_2\ge 0$.

Our main theorem is as follows:

\begin{theorem}\label{thm:main1}
Let the assumptions be given as above, and $\varphi$ is a $C^2$ solution to the equation \eqref{eqn:main}, then there exists a constant $C>0$ depending on the given parameters $n, A, K, p, \gamma$, and $\delta, R, \kappa_1, \kappa_2$ such that 
\begin{equation}\label{eqn:main1}
{\mathrm{sup}}_X | (\varphi + \theta) - {\mathrm{sup}}_X(\varphi + \theta)   | \le C.
\end{equation}

\end{theorem}

To prove this theorem, we need the following preliminary lemmas. First, 
we recall an estimate on the Green's function $G_\omega$ associated to $\omega$.

\begin{lemma}[\cite{GPSS22, GPS}]
\label{lemma 3.0}
There exists a constant $C_0 >0$ which depends on $p, A, K,\gamma$ and $\omega_X$ such that 
\begin{equation}\label{eqn:Green}
\sup_{x\in X} \Big( \sup_{y\in X}  (-G_\omega(x,y)) \cdot V_\omega + \int_X |G_\omega(x,y)| \omega^n(y)   \Big) \le C_0.
\end{equation}  
\end{lemma}

Next we have the following uniform $L^1$ and $L^q$ bounds:

\begin{lemma}
\label{lemma 3.01}
Let $\omega\in {\mathcal W}$ be a K\"ahler metric, then we have
\begin{equation}\label{eqn:L1}
\frac{1}{V_\omega}\int_X (-\psi) \omega^n \le 2n C_0,\quad \forall \psi\in  PSH(X,\omega)\mbox{ with }\sup_X \psi = 0,
\end{equation}
where $C_0$ is the constant in Lemma \ref{lemma 3.0}.
\end{lemma}
\noindent{\em Proof of Lemma \ref{lemma 3.01}.} By smooth approximation we may assume $\psi$ is $C^2$ and $\omega$-PSH. We assume $\sup_X \psi = \psi(x_0) = 0$ for some $x_0\in X$. Applying Green's formula to $\psi$ at $x_0$ gives
$$\psi(x_0) = \frac{1}{V_\omega} \int_X \psi \omega^n + \int_X [G_\omega(x_0,\cdot) - \inf_{z\in X} G_\omega(x_0, z)] (-\Delta _\omega \psi) \omega^n.$$
Hence by Lemma \ref{lemma 3.0}, we have
$$\frac{1}{V_\omega} \int_X (- \psi) \omega^n \le n \int_X \{|G_\omega(x_0,\cdot)| + \frac{C_0}{V_\omega} \}\omega^n \le 2n C_0.$$
\hfill Q.E.D.

\begin{lemma}\label{lemma 3.1}
Let $\omega\in {\mathcal W}$ be a K\"ahler metric. For any $q>n$, there exists a constant  $C_1 = C_1(p, A, K,\gamma, \omega_X,q)>0$ such that for any  $\omega$-PSH function $\psi$ with $\sup_X \psi = 0$, we have
\begin{equation}\label{eqn:alpha type}
\frac{1}{V_\omega}\int_X ({- \psi})^q \omega^n \le C_1. 
\end{equation}
\end{lemma}
\noindent {\em Proof of Lemma \ref{lemma 3.1}.} By Proposition 3.1 in \cite{GPSS22}, the condition $[\omega]\cdot[\omega_X]^{n-1}\le A$ implies that there exists a smooth $(1,1)$-form $\tilde \omega\in [\omega]$ such that \begin{equation}\label{eqn:lemma 0a}- C_2 \omega_X \le \tilde\omega \le C_2 \omega_X\end{equation} for some $C_2 = C_2(A,\omega_X)>0$. By $\ddbar $-lemma, we can write $\omega = \tilde \omega + \ddbar u$ for an $\tilde \omega$-PSH function $u$ normalized by $\sup_X u = 0$. 
%
Since by \eqref{eqn:lemma 0a} it holds that $\tilde \omega\le C_2 \omega_X$, we see that $u$ is  $(C_2 \omega_X)$-PSH, hence the $\alpha$-invariant estimates hold for $u$ (see \cite{H, T}), that is, there exist $\alpha' = \alpha'(C_2, \omega_X)>0$ and $C_3 = C_3(C_2, \omega_X)>0$ such that 
\begin{equation}\label{eqn:lemma a}
\int_X e^{-\alpha' u} \omega_X^n \le C_3.
\end{equation} 
For any $\omega$-PSH function $\psi$ with $\sup_X \psi = 0$, $\psi + u$ is also $\tilde \omega$-PSH, hence $(C_2\omega_X)$-PSH. Again by $\alpha$-invariant estimates, we have
\begin{equation}\label{eqn:lemma b}
\int_X e^{\alpha'( \sup_X(\psi + u) - \psi -  u)} \omega_X^n \le C_3.
\end{equation} 
To finish the proof, it suffices to prove a lower bound of $\sup_X (\psi + u)$. We note that by Young's inequality
\begin{equation}\label{eqn:lemma c}
\frac{1}{V_\omega}\int_X (-u)  \omega^n = \int_X (-u) e^{F_\omega} \omega_X^n \le C \int_X e^{F_\omega} (1+ |F_\omega|) \omega_X^n + C \int_X e^{-\alpha' u} \omega_X^n\le C_4,
\end{equation}
where we have used \eqref{eqn:lemma a} and $C_4>0$ depends on $p, A, K,$ and $\omega_X$. Combining \eqref{eqn:L1} and \eqref{eqn:lemma c} yields $$
\frac{1}{V_\omega}\int_X (-u - \varphi)  \omega^n \le 2nC_0 + C_4 = C_5.
$$
Hence $\sup_X(u+\varphi) \ge - C_5$ as claimed. Together with \eqref{eqn:lemma b} this implies that
\begin{equation}\label{eqn:lemma d}
\int_X e^{\alpha'(  - \psi )} \omega_X^n \le C_3 e^{C_5} = : C_6.
\end{equation} To prove \eqref{eqn:alpha type}, we apply Young's inequality again to get
$$(-\psi) e^{F_\omega/q} \le e^{F_\omega/q} ( 1 + |F_\omega|^{p/q}) + C_{p,q} e^{-\alpha' \psi/q}.$$
Taking $q$th power on both sides of the above inequality gives
\begin{equation}\label{eqn:lemma e}
(-\psi)^q e^{F_\omega} \le C e^{F_\omega} ( 1 + |F_\omega|^{p}) + C_{p,q} e^{-\alpha' \psi}.
\end{equation}
Integrating both sides of \eqref{eqn:lemma e} over $X$ against the measure $\omega_X^n$ yields the desired inequality \eqref{eqn:alpha type}, by making use of \eqref{eqn:lemma d}. \hfill Q.E.D.

\medskip



\begin{lemma}\label{lemma 3.2}
Let the assumptions be as above. For any $v\in C^2$ with $\sup_X v = 0$ satisfying $\lambda[h_{\chi' + \ddbar v}  ]\in \Gamma\subset \Gamma_1$, we have 
\begin{equation}\label{eqn:lemma 3.2}
\frac{1}{V_\omega} \int_X (-v) \omega^n \le 2 \kappa_2 C_0 =: C_7,
\end{equation}
where $C_0$ is the constant in Lemma \ref{lemma 3.0}.

\end{lemma}
\noindent{\em Proof of Lemma \ref{lemma 3.2}.} The condition $\lambda[h_{\chi' + \ddbar v}  ]\in \Gamma_1$ implies that
$\tr_\omega \chi' + \Delta_\omega v >0$. By \eqref{eqn:chi} this gives $\Delta_\omega v\ge -\kappa_2$. Then we can apply the Green's formula for $v$ at its maximum point $x_0$ to get
$$0 \le \frac{1}{V_\omega} \int_X v \omega^n + \int_X [ G_\omega(x_0,\cdot) - \inf_z G_\omega(x_0,z)   ] \kappa_2\omega^n. $$
\hfill Q.E.D.


\noindent {\em Proof of Theorem \ref{thm:main1}.}  Let $\varphi$ be a $C^2$ function solving the equation \eqref{eqn:main}. Denote $\varphi ' = \varphi + \theta$ and by adding an appropriate constant to $\varphi'$ we may assume $\sup_X\varphi ' = 0$. We have
$$\chi + \ddbar \varphi = \chi' + \ddbar \varphi'.$$
Hence, we can view $\varphi'$ as a solution to \eqref{eqn:main} if $\chi$ is replaced by $\chi'$.
As in \cite{GPT, GP22a}, for any $s\ge 0$, we consider the sublevel set $\Omega_s$ of $\varphi'$ by $\Omega_s = \{z\in X: \varphi'(z) < -s\}$. We solve the following auxiliary complex Monge-Amp\`ere equation which admits a smooth solution since $\omega$ is a K\"ahler metric \cite{Y} 
\begin{equation}\label{eqn:aux}
(\omega + \ddbar \psi_{s,k})^n = \frac{\tau_k(-\varphi '- s)}{A_{s,k}}  \omega^n, \quad \sup_X \psi_{s,k} = 0,
\end{equation}
where the constant $A_{s,k}$ is defined so that the equation is compatible, i.e. 
\begin{equation}\label{eqn:Ask}
A_{s,k}: = \frac{1}{V_\omega} \int_X \tau_k(-\varphi '- s)  \omega^n \to \frac{1}{V_\omega} \int_{\Omega_s} (-\varphi '- s)  \omega^n=: A_s,
\end{equation}
as $k\to\infty$. We consider the comparison function 
\begin{equation}\label{eqn:test fn}\Psi = -\varepsilon (-\psi_{s,k} + \Lambda)^{\frac{n}{n+1}} - \varphi '- s,\end{equation}
where \begin{equation}\label{eqn:epsilon}
\varepsilon = \beta_1 A_{s,k}^{\frac{1}{n+1}}, \quad \mbox{where } \beta_1:=\Big(\frac{n+1}{n} (R+\delta +\kappa_1)\Big)^{\frac{n}{n+1}},\end{equation}
and $\Lambda$ is given by 
\begin{equation}\label{eqn:Lambda} \Lambda =\beta_2 A_{s,k},\quad \mbox{where }\beta_2 = \frac{n}{n+1} (\frac{1}{\delta})^{n+1} (R+\delta + \kappa_1)^n.\end{equation}
We {\bf claim} that $\Psi\le 0$ on $X$. By definition, it suffices to check this when the maximum point of $\Psi$ is achieved at $x_0\in \Omega_s$. By the maximum principle, we have $\ddbar\Psi\le 0$ at $x_0$, which subsequently implies that at $x_0$,
$$\ddbar \varphi '\ge \frac{n\varepsilon}{n+1} (-\psi_{s,k} + \Lambda)^{-\frac{1}{n+1}} \ddbar \psi_{s,k} + \frac{n\varepsilon}{(n+1)^2} (-\psi_{s,k} + \Lambda)^{-\frac{n+2}{n+1}} \partial \psi_{s,k}\wedge \bar \partial \psi_{s,k}. $$
Therefore, adding $\chi'$ on both sides of the above we obtain
\bea \nonumber
\chi' + \ddbar \varphi' &\ge & \frac{n\varepsilon}{n+1} (-\psi_{s,k} + \Lambda)^{-\frac{1}{n+1}} (\omega + \ddbar \psi_{s,k}) \\
&& \nonumber + \chi' +  [ -  \frac{n\varepsilon}{n+1} (-\psi_{s,k} + \Lambda)^{-\frac{1}{n+1}}  ]\omega\\
&\ge & \label{eqn:key} \frac{n\varepsilon}{n+1} (-\psi_{s,k} + \Lambda)^{-\frac{1}{n+1}} (\omega+ \ddbar \psi_{s,k}) +\chi' - \delta \omega, 
\eea
where in the last line we have applied the choice of constant $\Lambda$ in \eqref{eqn:Lambda} and the assumption \eqref{eqn:chi}. From  the inequality \eqref{eqn:key}, we readily get that at the point $x_0$, 
\begin{equation}
\label{eqn:key 1}
\chi ' + \ddbar\varphi ' = \chi' - \delta \omega + \rho,
\end{equation}
for some {\em positive} $(1,1)$-form $\rho$. In terms of eigenvalues, \eqref{eqn:key 1} is equivalent to 
$$\lambda[h_{\chi' + \ddbar \varphi'}] \in \lambda[h_{\chi'}] - \delta {\mathbf 1} + \Gamma_n.$$ Moreover, from the equation \eqref{eqn:main}  we also have$f(\lambda[h_{\chi' + \ddbar \varphi'}] ) = e^{F(x_0)}$. Thus by the definition of a $(\delta,R)$-subsolution  $\chi'$, we conclude that the vector of eigenvalues $\lambda[h_{\chi' + \ddbar \varphi'}]|_{x_0} \in B(0, R)\subset {\mathbb R}^n$. If we denote $$\lambda[h_{\chi' + \ddbar \varphi'}]  |_{x_0}= (\lambda_1,\ldots,\lambda_n),$$ then $|\lambda_i|\le R$ for each $i = 1,\ldots, n$. Taking trace with respect to $\omega$ on both sides of \eqref{eqn:key},  together with the assumption \eqref{eqn:chi}, we obtain that
\begin{equation}\label{eqn:key 2}
n R + n\delta + n \kappa_1 \ge  \frac{n\varepsilon}{n+1} (-\psi_{s,k} + \Lambda)^{-\frac{1}{n+1}} \tr_\omega(\omega + \ddbar \psi_{s,k})>0.
\end{equation}
Applying the arithmetic-geometric inequality to the RHS of \eqref{eqn:key 2}, and using the auxiliary equation \eqref{eqn:aux} satisfied by $\omega + \ddbar \psi_{s,k}$, we get
\bea
\nonumber
(nR+n \delta + n\kappa_1) &\ge & \frac{n^2\varepsilon}{n+1} (-\psi_{s,k} + \Lambda)^{-\frac{1}{n+1}} \Big(\frac{\tau_k(-\varphi' - s)}{A_{s,k}} \Big)^{1/n}\\
\label{eqn:key 3} &\ge &\frac{n^2\varepsilon}{n+1} (-\psi_{s,k} + \Lambda)^{-\frac{1}{n+1}} \frac{(-\varphi' - s)^{1/n}}{A^{1/n}_{s,k}}.
\eea
From the choice of the constant $\varepsilon$ in \eqref{eqn:epsilon}, the inequality \eqref{eqn:key 3} easily yields that $\Psi(x_0)\le 0$. Hence we finish the proof the claim. We then obtain that on $\Omega_s$,
\begin{equation}\label{eqn:key 4}{(- \varphi' - s)^{\frac{n+1}{n}}}\le \beta_1^{\frac{n+1}{n}} {A_{s,k}^{1/n}} (-\psi_{s,k} + \Lambda). \end{equation}
We fix a constant $q>n$, and by Lemma \ref{lemma 3.1}, there exists a constant $C_1>0$ such that $\frac{1}{V_\omega}\int_X (-\psi_{s,k})^q \omega^n\le C_1$. Taking $q$-th power on both sides of \eqref{eqn:key 4} and integrating the resulted inequality over $\Omega_s$ yield that 
\bea\label{eqn:key 6}\nonumber
\frac{1}{V_\omega}\int_{\Omega_s}  (- \varphi' - s)^{q\frac{n+1}{n}} \omega^n & \le & C_2 {A_{s,k}^{q/n}} \frac{1}{V_\omega} \int_X ((-\psi_{s,k})^q + \Lambda^q) \omega^n\\
\label{eqn:key 6} &\le & C_2(C_1 + \Lambda^q) A_{s,k}^{q/n},
\eea
where $C_2$ is a constant depending on $q$ and $\beta_1$. Letting $k\to \infty$ in \eqref{eqn:key 6} we see that this inequality remains true if $A_{s,k}$ is replaced by $A_s$ (along with the $A_{s,k}$ in the constant $\Lambda$ in \eqref{eqn:Lambda}).

 We define a monotone function as in \cite{GPT} $$\phi: {\mathbb R}_{\ge 0}\to {\mathbb R}_{\ge 0}\quad \mbox {by }\quad \phi(s) = \frac{1}{V_\omega} \int_{\Omega_s} \omega^n.$$
From H\"older inequality, it follows that
\bea\label{eqn:H need}
A_s & = &\frac{1}{V_\omega}\int_{\Omega_s} (-\varphi' - s) \omega^n \\
 &\le & \nonumber \Big(\frac{1}{V_\omega} \int_{\Omega_s}  (-\varphi' -s)^{\frac{q(n+1)}{n}}\omega^n \Big)^{\frac{n}{q(n+1)}} \phi(s)^{\frac{q(n+1) - n}{q(n+1)}}\\
&\le & (C_2(C_1 + \Lambda^q) A_{s}^{q/n}) ^{\frac{n}{q(n+1)}} \phi(s)^{\frac{q(n+1)}{q(n+1) - n}} \nonumber \\  &=& \nonumber (C_2(C_1 + \Lambda^q)) ^{\frac{n}{q(n+1)}}A_{s}^{\frac{1}{n+1}} \phi(s)^{\frac{q(n+1) - n}{q(n+1)}}.
\eea
Hence we have 
\begin{equation}\label{eqn:key 7}
A_{s}\le (C_2(C_1 + \Lambda^q)) ^{\frac{1}{q}} \phi(s)^{1 + \frac{1}{n} - \frac{1}{q}} = C_3 \phi(s)^{1 + \frac{1}{n} - \frac{1}{q}},
\end{equation} where the constant $C_3$ is given by $ (C_2(C_1 + \Lambda^q)) ^{\frac{1}{q}}$ and is uniformly bounded since $A_s$ is bounded from Lemma \ref{lemma 3.2}. 
From \eqref{eqn:key 7} we get 
\begin{equation}\label{eqn:key 8}
 r\phi(s+r) \le C_3  \phi(s)^{1 + \frac{1}{n} - \frac{1}{q}},\quad \forall r\in (0,1) \mbox{ and }s\ge 0.
\end{equation}
It is a classic fact due to De Giorgi \cite{K}, that \eqref{eqn:key 8} implies $\Omega_{S_\infty} = \emptyset$, for $S_\infty = \frac{1}{1 - 2^{-nq/(q-n)}} + s_0,$  where $s_0>0$ is a number that satisfies $\phi(s_0)\le (2C_3)^{(q-n)/nq}$. 

We note that 
$$\phi(s_0) \le \frac{1}{s_0} \frac{1}{[\omega]^n} \int_X (-\varphi') \omega^n \le  \frac{C_7}{s_0}, $$ where the last inequality follows from Lemma \ref{lemma 3.2} and $C_7>0$ is the constant in that lemma. So we can take $s_0 =C_7 (2C_3)^{ - (q-n)/nq}$. Combining the above, we conclude that 
\begin{equation}\label{eqn:key 9}
\sup_X (-\varphi')\le S_\infty = C_7 (2C_3)^{ - (q-n)/nq}  + \frac{1}{1 - 2^{-nq/(q-n)}}. 
\end{equation}
This completes the proof of Theorem \ref{thm:main1}. \hfill Q.E.D.

\medskip

\noindent{\bf Remark 3.1.} We only require the K\"ahler metric $\omega$ to be in the set ${\mathcal W}$ for the purpose of verifying the two integrability conditions stated in Lemma \ref{lemma 3.1} and Lemma \ref{lemma 3.2}. It is important to note that the proof of Theorem \ref{thm:main1} remains valid when these two lemmas are satisfied. Therefore, there exist alternative choices for the family of metrics $\omega$ beyond the set $\mathcal W$ that can ensure the validity of these lemmas, and consequently, Theorem \ref{thm:main1}. For instance, in case $\Gamma = \Gamma_n$, we can impose the following assumptions on $\omega$ and the subsolutions $\chi'$: 
\begin{equation}\label{eqn:extra}\omega\le K_1 \omega_X,\quad \frac{1}{V_\omega}\frac{\omega^n}{\omega^n_X}\le K_2, \quad -\kappa_1 \omega \le \chi'\le \kappa_2 \omega.\end{equation}
Under these assumptions, Lemma \ref{lemma 3.1} holds from the $\alpha$-invariant estimates readily. Note that $\chi' + \ddbar \varphi' >0$ and hence $\kappa_2 K_2 \omega_X + \ddbar \varphi' >0$. Green's formula for the fixed metric $\omega_X$ shows $\int_X (-\varphi') \omega_X^n$ is bounded depending only on the given parameters. Lemma \ref{lemma 3.2} is also true by the second inequality in \eqref{eqn:extra}.

\section{The Hermitian case}
\setcounter{equation}{0}
In this section, we will prove the $L^\infty$ estimates when the metric $\omega$ is {\em fixed} and only Hermitian. Let $X$ be a compact complex manifold without boundary. As in the last section, we assume that $\chi' = \chi - \ddbar \theta$ is a $(\delta, R)$-subsolution to \eqref{eqn:main} which satisfies 
\begin{equation}\label{eqn:new chi}
-n \kappa_1 \omega \le \chi' \le \kappa_2 \omega.
\end{equation}

\begin{theorem}\label{thm:main2}
Let the assumptions be given as above, and $\varphi$ is a $C^2$ solution to the equation \eqref{eqn:main}, then there exists a constant $C>0$ depending only on $n, \omega$, and $\delta, R, \kappa_1, \kappa_2$ such that 
\begin{equation}\label{eqn:main1}
{\mathrm{sup}}_X | (\varphi + \theta) - {\mathrm{sup}}_X(\varphi + \theta)   | \le C.
\end{equation}

\end{theorem}
We denote $\varphi ' = \varphi +\theta$ and assume without loss of generality that $\sup_X \varphi' = 0$. 
Since the metric $\omega$ is fixed, there is a constant $r_0 = r_0(n,\omega, X)>0$ such that for any point $z_0\in X$, there exist holomorphic coordinates $z = (z_1,\ldots, z_n)$ at $z_0$ such that on the ball $B(z_0, 2r_0) = \{z:~ \sum_j |z_j|^2 < (2r_0)^2\}$ the metric $\omega = g_{i\bar j} \sqrt{-1} dz_i\wedge d\bar z_j$ satisfies 
\begin{equation}\label{eqn:omega 1}
\frac{1}{2} \delta_{ij} \le g_{i\bar j} \le 2 \delta_{ij}.
\end{equation}
We assume $x_0$ is a minimum point of $\varphi'$, i.e. $\min _X \varphi '= \varphi'(x_0)$, and pick the coordinates $z$ at $x_0$ such that \eqref{eqn:omega 1} holds. Fix $s_0 = \frac{\delta}{10} \kappa_1 r_0^2>0$. For any $s\in (0, s_0)$, we consider the function $u_s$ defined on the ball $B(x_0, r_0) = \{z~: ~|z|<r_0\}$ as
$$u_s(z) = -\varphi'(x_0) + \varphi'(z) + \frac{\delta}{4}  |z|^2 -  s,$$
 and define the level set of $u_s$ 
 $$U_s: = \{z\in B(x_0, r_0)| ~ u_s(z) < 0\}.$$
 Since $u_s|_{\partial B(x_0, r_0)}>0$, we see that $U_s$ is relatively compact in $B(x_0, r_0)$. We consider the solution to the Dirichlet problem
$$(\ddbar \psi_{s,k})^n = \frac{\tau_k(-u_s)}{A_{s,k}}  \omega^n,\, \mbox{in }B(x_0, r_0),$$
with $\psi_{s,k} = 0$ on $\partial B(x_0,r_0)$, where we choose the constant $A_{s,k}>0$ such that \begin{equation}\label{eqn:Diri}\int_{B(x_0, r_0)} (\ddbar \psi_{s,k})^n = 1\end{equation} in other words 
$$A_{s,k}: = \int_{B(x_0, r_0)} \tau_k(-u_s) \omega^n \to A_s := \int_{\Omega_s} (-u_s) \omega^n.$$
Note that the solution $\psi_{s,k}\le 0$ exists and is unique by \cite{CKNS}. We next consider the test function
$$\Phi = -\varepsilon (-\psi_{s,k})^{\frac{n}{n+1}} - u_s, \mbox{ in }B(x_0, r_0),$$
where \begin{equation}\label{eqn:new epsilon}\varepsilon = (\frac{n+1}{n})^{\frac{n}{n+1}} (R+\delta + \kappa_1)^{\frac{n}{n+1}} A_{s,k}^{\frac{1}{n+1}} =:\beta_3 A_{s,k}^{\frac{1}{n+1}}.\end{equation}
As in the last section, we claim $\Phi\le 0$. It suffices to prove this when the maximum of $\Phi$ is achieved at some point $z_0\in U_s\subset B(x_0, r_0)$. By maximum principle, we have that   at the point $z_0$, $\ddbar \Phi \le 0$, which implies that at $z_0$
\bea \nonumber
\ddbar \varphi '& \ge& \frac{n\varepsilon}{n+1} (-\psi_{s,k})^{-\frac{1}{n+1}} \ddbar \psi_{s,k} - \frac{\delta}{4} \omega_{{\mathbf C}^n}\\
&\ge \nonumber &\frac{n\varepsilon}{n+1} (-\psi_{s,k})^{-\frac{1}{n+1}} \ddbar \psi_{s,k} - \frac{\delta}{2} \omega. 
\eea
Hence at $z_0$ it holds that
\begin{equation}\label{eqn:Hnew 1}
\chi' + \ddbar \varphi' > \frac{n\varepsilon}{n+1} (-\psi_{s,k})^{-\frac{1}{n+1}} \ddbar \psi_{s,k} + \chi' - {\delta} \omega. 
\end{equation}
This inequality shows that $\chi' + \ddbar \varphi' = \chi' - \delta \omega + \rho$ for some positive definite $(1,1)$-form $\rho$. In terms of eigenvalues, this means that $$\lambda[h_{\chi'+ \ddbar \varphi'}] \in \lambda[h_{\chi'} ] - \delta {\mathbf 1} + \Gamma_n.$$ By the definition of the $(\delta, R)$-subsolution $\chi'$, we have $|\lambda_i|\le R$ for each $i=1,2,\ldots, n$, where $(\lambda_1,\ldots, \lambda_n)$ are the eigenvalues of $\omega^{-1}\cdot(\chi' + \ddbar \varphi' )|_{z_0}$. Again by \eqref{eqn:Hnew 1} we have at $z_0$
\bea\nonumber 
n( R + \kappa_1 +  \delta)& \ge &  \tr_\omega( \chi' + \ddbar \varphi' -\chi' + \delta \omega   )\\
 & \ge & \nonumber \frac{n\varepsilon}{n+1} (-\psi_{s,k})^{-\frac{1}{n+1}} \tr_\omega (\ddbar \psi_{s,k} )\\
 &\ge &\nonumber \frac{n^2\varepsilon}{n+1} (-\psi_{s,k})^{-\frac{1}{n+1}} \Big (\frac{(\ddbar \psi_{s,k})^n}{\omega^n} \Big)^{1/n}\\
  &\ge &\nonumber \frac{n^2\varepsilon}{n+1} (-\psi_{s,k})^{-\frac{1}{n+1}} \frac{(-u_s)^{1/n}}{A_{s,k}^{1/n}},
\eea
which is equivalent to $\Phi|_{z_0}\le 0$ by the choice of $\varepsilon$ in \eqref{eqn:new epsilon}. Therefore, on $U_s$ we have 
\begin{equation}\label{eqn:omega 5}
{(-u_s)^{{(n+1)/}{n}}}{} \le \beta_3^{\frac{n+1}{n}}  A_{s,k}^{1/n} (-\psi_{s,k}).
\end{equation}
By the uniform $\alpha$-invariant estimates (see \cite{K}) for $\psi_{s,k}$ which satisfies \eqref{eqn:Diri}, there exist uniform constants $\alpha = \alpha(n,r_0, \omega)$ and $C= C(n,r_0, \omega)$ such that $\int_{B(x_0,r_0)} e^{-\alpha \psi_{s,k}}\omega^n \le C$, hence for any $q>n$ it holds that
\begin{equation}\label{eqn:Hnew 2}
\int_{B(x_0,r_0)} (-\psi_{s,k})^q\omega^n \le C_0' =  C_0'(n,r_0, \omega, q).
\end{equation}
Taking $q$-th power on both sides of \eqref{eqn:omega 5} and integrating over $U_s$, we have from \eqref{eqn:Hnew 2} that
\begin{equation}\label{eqn:Hnew 3}
\int_{U_s} (-u_s)^{\frac{q(n+1)}{n}} \omega^n \le \beta_3^{\frac{q(n+1)}{n}} A_{s,k}^{\frac{q}{n}} \int_{B(x_0,r_0)} (-\psi_{s,k})^q \omega^n \le C_0'\beta_3^{\frac{q(n+1)}{n}} A_{s,k}^{\frac{q}{n}}.
\end{equation}
Letting $k\to\infty$, inequality \eqref{eqn:Hnew 3} continues to hold when $A_{s,k}$ is replaced by $A_s$. 
By H\"older inequality, we can argue similarly as in \eqref{eqn:H need} to derive
\begin{equation}\label{eqn:Hnew 4}
A_{s}\le (C_0'\beta_3^{\frac{q(n+1)}{n}}) ^{\frac{1}{q}}\hat \phi(s)^{1 + \frac{1}{n} - \frac{1}{q}} =: C_1' \hat\phi(s)^{1 + \frac{1}{n} - \frac{1}{q}},
\end{equation}
where $C_1' =  (C_0'\beta_3^{\frac{q(n+1)}{n}}) ^{\frac{1}{q}}$ and $\hat\phi(s)$ is a monotone function given by
$$\hat \phi(s) = \int_{U_s} \omega^n,\quad \forall s\in (0,s_0).$$
From \eqref{eqn:Hnew 4} we obtain
\begin{equation}\label{eqn:Hnew 5}
t \hat\phi(s-t) \le C_1' \hat\phi(s)^{1+\frac{1}{n} - \frac{1}{q}},\quad \forall \,0 < t < s <s_0. 
\end{equation}
Then an iteration argument (see e.g. Lemma 7 in \cite{GP22a}) shows that there exists a constant $c_0>0$ that depends on $n, q, C_0',$ and $s_0$ such that $\hat\phi(s_0)\ge c_0$. We observe that on $U_{s_0}$, 
$$u_{s_0}(z)< 0 \quad \Rightarrow \quad -\varphi'(x_0) \le -\varphi'(z) - \frac{\delta}{4 } |z|^2 + s_0. $$
Integrating this inequality over $U_{s_0}$ yields 
\begin{equation}\label{eqn:Hnew 6}
c_0(-\varphi'(x_0)) \le \int_{U_{s_0}} (-\varphi'(x_0)) \omega^n\le \int_{B(x_0, r_0)} (-\varphi') \omega^n + C(n, r_0) s_0. 
\end{equation}
From \eqref{eqn:Hnew 6}, we see that
$$-\varphi'(x_0) = -{\mathrm{inf}}_X \varphi' \le C( \| \varphi'\|_{L^1(X,\omega^n)} + 1  ),$$
for some constant $C>0$ that depends on $n, \kappa_1, q$ and $\omega$. To complete the proof it suffices to derive an $L^1$ estimate of $\varphi'$. This follows from the fact that $\lambda[h_{\chi'+ \ddbar \varphi'}]\in \Gamma\subset \Gamma_1$, hence 
$$\tr_\omega \chi' + \Delta_{\omega} \varphi' \ge 0,$$ where $\Delta_\omega \varphi' = \tr_{\omega}(\ddbar \varphi')$ is the rough Laplacian of the Hermtian metric $\omega$. Given the assumption in \eqref{eqn:new chi} that $\tr_\omega \chi' \le \kappa_2$, we have the linear differential inequality $\Delta_{\omega} \varphi' \ge - \kappa_2$ on $X$. By a local covering argument (see e.g. Lemma 8 in \cite{GP22a}) and standard elliptic estimates \cite{GT}, we have $$\int_X (-\varphi')\omega^n \le C,$$
for some $C>0$ depending on $n,\omega, \kappa_2$. This finishes the proof of Theorem \ref{thm:main2}.

\section{Other forms of subsolutions}
In this section, we discuss another form of subsolutions in the case of $\Gamma = \Gamma_n$. To simplify the exposition, we make the convention that for a vector $\lambda = (\lambda_1,\ldots, \lambda_n)\in \Gamma_n$, i.e. $\lambda_i>0$ for all $i$, we define $\tilde \lambda = (\frac{1}{\lambda_1},\ldots, \frac{1}{\lambda_n})\in \Gamma_n$. Abusing notation, we sometimes write $\tilde \lambda_i = \frac{1}{\lambda_i}$. Associated with the operator $f$ which satisfies  the condition (ii) in Section 1 and the condition of homogeneity of degree one, we define a new fully nonlinear operator $\tilde f : \Gamma_n \to {\mathbf R}_{>0}$ as
\begin{equation}\label{eqn:tilde f}
\tilde f(\tilde \lambda): = \frac{1}{f({\tilde \lambda_1}^{-1},\ldots, {\tilde \lambda_n}^{-1})} = \frac{1}{f(\lambda)}.
\end{equation}
Clearly $\tilde f$ is also a homogeneous function of degree one and $\frac{\partial \tilde f}{\partial \tilde \lambda_i}(\tilde\lambda) >0$ for any $\tilde \lambda \in \Gamma_n$. Furthermore, we impose the extra conditions on $f$ or $\tilde f$ as follows.

\smallskip

(iii): $\tilde f$ is {\em concave} in $\tilde \lambda \in\Gamma_n$.

\begin{definition}\label{defn 5.1}
We say a metric $\chi' = \chi  - \ddbar \theta$ is a $\delta$-subsolution of (\ref{eqn:main}), if 
\begin{equation}\label{eqn:sub}
\sum_{j\neq i}\tilde \mu_j \frac{\partial \tilde f}{\partial \tilde \lambda_j}(\tilde\mu) \le (1-\delta) e^{-F(z)},\quad \forall z\in X, \, \forall i = 1,\ldots, n.
\end{equation}
%
Here $\mu = (\mu_1,\ldots, \mu_n)$ is an unordered vector of eigenvalues of $\omega^{-1}\cdot \chi'$, and conventionally we have defined $\tilde \mu_j  = \frac{1}{\mu_j}$ and $\tilde \mu  = (\tilde \mu_1,\ldots, \tilde \mu_n)$. 
\end{definition}
We observe that if $\chi'$ satisfies the equation (\ref{eqn:main}), by the homogeneity of degree one condition, we necessarily have $\sum_{j=1}^n\tilde \mu_j \frac{\partial \tilde f}{\partial \tilde \lambda_j}(\tilde\mu) = e^{-F(z)}$. So a solution is indeed a $\delta$-subsolution for some $\delta>0$.  

Given a $\delta$-subsolution $\chi'$ as in Definition \ref{defn 5.1}, we have the following uniform $L^\infty$ estimates for solutions $\varphi$ to \eqref{eqn:main}.

\begin{theorem}\label{lemma 5.1}

Assume the operator $f$ satisfies the additional condition (iii). Suppose $\varphi$ is a $C^2$ solution to the equation \eqref{eqn:main}, $\omega\in {\mathcal W}$ is a K\"ahler metric, and $\chi'$ is a $\delta$-subsolution as in Definition \ref{defn 5.1} such that $\chi'\le \kappa_2 \omega$ for some constant $\kappa_2>0$. Then there exists a constant $C>0$ that depends on $\kappa_2,\delta$ and the parameters in ${\mathcal W}$ such that
$${\mathrm{sup}}_X | (\varphi + \theta) -  {\mathrm{sup}}_X (\varphi + \theta)  |\le C.$$ 

\end{theorem}
\noindent {\em Proof of Theorem \ref{lemma 5.1}.} We follow the notations as in Section \ref{section Kahler}, and denote $\varphi ' = \varphi + \theta$. Let $\psi_{s,k}$ and $\Psi$ be as in \eqref{eqn:aux} and \eqref{eqn:test fn}, respectively. Here the two constants $\varphi$ and $\Lambda$ in \eqref{eqn:test fn} are chosen as
$$\varepsilon = \Big( \frac{n+1}{n} R \Big)^{\frac{n}{n+1}} A_{s,k}^{\frac{1}{n+1}},\quad \Lambda = \Big( \frac{2n\kappa_2 \varepsilon}{(n+1)\delta}\Big)^{n+1},  $$
and here $ R:= \frac{(2-\delta)(1-\delta)}{\delta}\kappa_2 >0$. With these choice of constants, we claim $\Psi\le 0$. At a maximum point $x_0\in\Omega_s$ of $\Psi$, an analogous inequality as \eqref{eqn:key} holds, which gives 
\begin{equation}\label{eqn:n c1}
\chi' + \ddbar \varphi' > (1- 0.5\delta) \chi'. 
\end{equation}
We take normal coordinates of $\omega$ at $x_0$, and consider  the Hermitian matrices $B' = (\chi'+ \ddbar \varphi')_{i\bar j}|_{x_0}$ and $B = (\chi ')_{i\bar j}|_{x_0}$. $B'$ and $B$ may not be necessarily diagonal under these coordinates. We may assume there eigenvalues are listed in an increasing order, i.e. 
$$0<\lambda_1\le \cdots \le \lambda_n, \quad\mbox{and}\quad  0< \mu_1\le \cdots \le \mu_n.$$ The inequality \eqref{eqn:n c1} means $B' >  (1-\delta)B.$ It follows from an elementary fact from linear algebra (see Lemma \ref{lemma 5.2} below) that 
\begin{equation}\label{eqn:n c2}
\lambda_i \ge  \mu_i(1- 0.5\delta),\quad \forall \, i = 1,\ldots, n.
\end{equation}
By the condition (i') that $\tilde f(\tilde \lambda)$ is concave in $\tilde \lambda$, it follows that
\begin{equation}\label{eqn:concave 1}
\tilde f(\tilde \lambda) \le \tilde f (\tilde \mu) + \sum_{i=1}^n \frac{\partial \tilde f}{\partial \tilde \lambda_i} (\tilde \mu) \cdot ( \tilde \lambda_i - \tilde \mu_i  ) =  \sum_{i=1}^n \frac{\partial \tilde f}{\partial \tilde \lambda_i} (\tilde \mu) \cdot  \tilde \lambda_i,
\end{equation}
where we have used the condition that $\tilde f(\tilde \lambda)$ is homogeneous of degree one. Here we need the subsolution condition \eqref{eqn:sub} for $\chi'$. For each fixed $i = 1,\ldots, n$, we have from \eqref{eqn:concave 1} that
\begin{equation}\label{eqn:concave 2}
e^{-F(x_0)} \le  \frac{\partial \tilde f}{\partial \tilde \lambda_i} (\tilde \bmu) \cdot  \tilde \lambda_i +   \sum_{j\neq i} \frac{\partial \tilde f}{\partial \tilde \lambda_j} (\tilde \bmu) \cdot  \tilde \lambda_j\le \frac{\partial \tilde f}{\partial \tilde \lambda_i} (\tilde \bmu) \cdot  \tilde \lambda_i +   \frac{1-\delta}{1-0.5\delta} e^{-F(x_0)},
\end{equation}
which implies that
$$\frac{\delta}{2 - \delta} e^{-F(x_0)}\le   \frac{\partial \tilde f}{\partial \tilde \lambda_i} (\tilde \bmu) \cdot  \tilde \lambda_i \le \frac{\tilde \lambda_i}{\tilde \mu_i} (1-\delta) e^{-F(x_0)}.$$ 
Hence $$\lambda_i \le \frac{(2-\delta)(1-\delta)}{\delta} \mu_i\le \frac{(2-\delta)(1-\delta)}{\delta}\kappa_2 =R.$$
The remaining arguments are the same as those in Section \ref{section Kahler}.  The proof of Theorem \ref{lemma 5.1} is completed. \hfill Q.E.D.

\bigskip

\begin{lemma}\label{lemma 5.2}
Suppose $A, B$ are both $n\times n$  Hermitian matrices satisfying $B\ge A$, i.e. $B-A$ is semi-positive definite. If we denote $\lambda_i(A)$ the eigenvalues of $A$ in the increasing order, i.e.  $ 0< \lambda_1(A)\le \cdots \le \lambda_n(A)$, and similar notation for $\lambda_i(B)$, then $$\lambda_i (B ) \ge \lambda_i(A),\quad \forall i = 1, \ldots, n.$$ 
\end{lemma}
\noindent{\em Proof of Lemma \ref{lemma 5.2}.} It follows from the Courant-Fischer min-max Theorem that 
$$\lambda_i(A) =\min\Big\{ \max \{ \langle Az, z \rangle: ~ z\in U,\, |z| = 1  \}, ~ U \mbox { with dim}_{{\mathbf C}} U = i \Big\},$$
where the minimum is taken over all {\em complex} vector spaces $U\subset {\mathbf C}^n$ with dimension $i$. For each fixed vector space $U$ we have 
$$\max \{ \langle Bz, z \rangle: ~ z\in U,\, |z| = 1  \} \ge \max \{ \langle Az, z \rangle: ~ z\in U,\, |z| = 1  \}.$$
Taking minimum over such $U$'s gives the desired inequality. \hfill Q.E.D.

\subsection{Examples}

\subsubsection{The $J$-equation} We consider the following form of $J$-equation:
\begin{equation}\label{eqn:J}
(\chi + \ddbar \varphi)^n = n e^{F} \omega\wedge (\chi + \ddbar \varphi)^{n-1}, \quad\sup_X \varphi = 0.
\end{equation} 
In the notation of $f(\bl[h_\varphi])$, this equation is equivalent to 
$$f(\bl) = \frac{\sigma_n (\bl)}{\sigma_{n-1} (\bl) } = e^F,$$
where $\bl\in \Gamma_n$ is the vector of eigenvalues of $\omega^{-1}\cdot (\chi+ \ddbar \varphi)$. Rewriting the nonlinear operator $f$ in terms of $\tilde f$, we see that $\tilde f(\tilde \bl) = \sigma_1(\tilde \bl)$, which satisfies the required conditions. The subsolution condition \eqref{eqn:sub} for $\chi'$ with this equation states that
$$\sum_{j\neq i} \frac{1}{\mu_j} \le (1-\delta){e}^{-F(z)},\quad \forall\, i = 1,\ldots, n,$$
which is just the subsolution condition introduced in \cite{SW} (see also \cite{Su}). We mention that for the $J$-equation, it is proved in \cite{So} that the existence of a subsolution is equivalent to a Nakai-Moishezon criterion, thus confirming a conjecture of \cite{LSz}. 

\subsubsection{Hessian type equations} The following equation is studied in \cite{FLM}, when $1\le k\le n$
$$
(\chi + \ddbar \varphi)^n = e^{F} \omega^{n-k}\wedge (\chi + \ddbar \varphi)^k,\quad \sup_X\varphi = 0.
$$
The subsolution condition \eqref{eqn:sub} for $\chi'$ with this equation is the one studied in \cite{FLM}.

\end{document}